\documentclass[11pt,twoside,reqno]{amsart}

\usepackage{amsmath}
\usepackage{amsthm}
\usepackage{amsfonts,amssymb}
\usepackage{bbm}
\usepackage{latexsym}
\usepackage{mathrsfs}
\usepackage[all]{xy}
\usepackage{url}
\usepackage[pdftex]{graphicx}
\usepackage{float}
\usepackage{tikz-cd}
\usepackage{hyperref}

\theoremstyle{plain}
\newtheorem{lema}{Lemma}[section]
\newtheorem{prop}[lema]{Proposition}
\newtheorem{teo}[lema]{Theorem}

\newtheorem{preg}[lema]{Question}
\newtheorem{coro}[lema]{Corollary}
\theoremstyle{remark}

\theoremstyle{definition}
\newtheorem{defi}[lema]{Definition}

\newcommand{\tr}{\mathrm{tr}} 
\newcommand{\rk}{\mathrm{rk}}

\newcommand{\Z}{\mathbb{Z}}

\newcommand{\PP}{\mathcal{P}}
\newcommand{\QQ}{\mathcal{Q}}

\newcommand{\id}{\mathbbm{1}}

\begin{document}

\title[Presentation complexes with the fixed point property]{Presentation complexes with the fixed point property}

\author[I. Sadofschi Costa]{Iv\'an Sadofschi Costa}

\address{Departamento de Matem\'atica\\
FCEyN-Universidad de Buenos Aires\\
Buenos Aires, Argentina}

\email{ivansadofschi@gmail.com}

\begin{abstract}
We prove that there exists a compact two-dimensional polyhedron with the fixed point property and even Euler characteristic. This answers a question posed by R.H. Bing in 1969. We also settle another of Bing's questions.
\end{abstract}



\maketitle

\section{Introduction}

In his influential article ``The elusive fixed point property" \cite{Bing}, R.H. Bing stated twelve questions. Since then eight of these questions have been answered \cite{Hagopian}. In this paper we answer Questions 1 and 8.

Recall that a space $X$ is said to have the \textit{fixed point property} if every map $f:X\to X$ has a fixed point. Motivated by an example of W. Lopez \cite{Lopez}, Bing stated in \cite{Bing} the following question.

\begin{preg}[Bing's Question 1]\label{Bing1} Is there a compact two-dimensional polyhedron with the fixed point property which has even Euler characteristic?
\end{preg}
This question was studied in \cite{Waggoner2}. In \cite{BarmakSadofschi} it was shown that such a space cannot have abelian fundamental group. In Corollary \ref{CoroBing1} we show that the answer to Question \ref{Bing1} is affirmative. Bing's Question 8 \cite{Bing} may be rephrased as follows.
\begin{preg}[Bing's Question 8]\label{Bing8}
What is the lowest dimension for a compact polyhedron $X$ with the fixed point property and such that a space $Y$ without the fixed point property can be obtained by attaching a disk $D$ to $X$ along an arc?
\end{preg}
The answer to this question is clearly greater than $1$. A one-dimensional polyhedron $X$ with the fixed point property is a tree, and then any space $Y$ obtained by attaching a disk along an arc is a contractible polyhedron.
According to C.L. Hagopian \cite{Hagopian}, Bing conjectured that the answer to Question \ref{Bing8} is $2$. This is the content of Theorem \ref{Bing8Example}.

\bigskip

\textbf{Acknowledgment:} I am grateful to Jonathan Barmak, without his advice and suggestions this paper would not have been possible.

\section{Bing groups}
If $\PP$ is a presentation, the presentation complex of $\PP$ will be denoted by $X_\PP$. Presentation complexes are in fact polyhedra. If a finite group $G$ is presented by a presentation $\PP$ with $g$ generators and $r$ relators, then $r-g$ is at least the number of invariant factors of $H_2(G)$. If this lower bound is attained for $\PP$, then the presentation is said to be \textit{efficient}.

\begin{defi} Let $G$ be a finite group and $d_1\mid \ldots \mid d_k$ be the invariant factors of $H_2(G)$. We say that $G$ is a \textit{Bing group} if for every endomorphism $\phi:G\to G$ we have $\tr(H_2(\phi)\otimes \id_{\Z_{d_1}})\neq -1$ in $\Z_{d_1}$.
\end{defi}

\begin{teo}\label{BingGroups} If $\PP$ is an efficient presentation of a Bing group $G$ then $X_{\PP}$ has the fixed point property.
\begin{proof}
Let $X=X_\PP$ and  $f:X\to X$ be a map. There is a $K(G,1)$ space $Y$ with $X=Y^2$. Now $f$ extends to a map  $\overline{f}:Y\to Y$. 
% We have a commutative diagram
% \begin{center}\begin{tikzcd}
% X\arrow[]{d}[swap]{f}\arrow[hook]{r}[]{i} & Y\arrow[]{d}[]{\overline{f}} \\
% X\arrow[hook]{r}[swap]{i} & Y
% \end{tikzcd}\end{center}
In the following commutative diagram, the horizontal arrows, induced by the inclusion $i:X\hookrightarrow Y$, are epimorphisms: 
\begin{center}\begin{tikzcd}
H_2(X)\arrow[]{d}[swap]{f_*}\arrow[two heads]{r}[]{i_*} & H_2(Y)\arrow[]{d}[]{\overline{f}_*} \\
H_2(X)\arrow[two heads]{r}[swap]{i_*} & H_2(Y)
\end{tikzcd}\end{center}
Let $d_1\mid\ldots \mid d_k$ be the invariant factors of $H_2(G)$. Since $\PP$ is efficient, the rank of $H_2(X)$ equals the number of invariant factors of $H_2(Y)$. Therefore the horizontal arrows in the following commutative diagram are isomorphisms:
\begin{center}\begin{tikzcd}
H_2(X)\otimes \Z_{d_1}\arrow[]{d}[swap]{f_*\otimes \id_{\Z_{d_1}}}\arrow[]{rr}[]{i_*\otimes \id_{\Z_{d_1}}}[swap]{\approx} && H_2(Y)\otimes \Z_{d_1}\arrow[]{d}[]{\overline{f}_*\otimes \id_{\Z_{d_1}}} \\
H_2(X)\otimes \Z_{d_1}\arrow[]{rr}[swap]{i_*\otimes \id_{\Z_{d_1}}}[]{\approx} && H_2(Y)\otimes \Z_{d_1}
\end{tikzcd}\end{center}

Now $\tr(f_*\otimes \id_{\Z_{d_1}})=\tr(\overline{f}_*\otimes \id_{\Z_{d_1}})\neq -1$ in $\Z_{d_1}$ since $G$ is a Bing group. Here we are using the natural isomorphism $H_2(BG)\approx H_2(G)$ of \cite[Theorem 5.1.27]{Rosenberg}. Recall that every map $BG\to BG$ is induced, up to homotopy, by an endomorphism $G\to G$.

Finally we obtain $\tr(f_*)\neq -1$ in $\Z$, since tensoring with $\Z_{d_1}$ reduces the trace modulo ${d_1}$.
So $L(f)\neq 0$ and, by the Lefschetz fixed point theorem, $f$ has a fixed point.
\end{proof}
\end{teo}

\begin{prop}\label{GIsBing} The group $G$ presented by $$\PP = \langle x,y \mid x^3,\, xyx^{-1}yxy^{-1}x^{-1}y^{-1},\, x^{-1}y^{-4}x^{-1}y^2x^{-1}y^{-1} \rangle$$ is a finite group of order $243$. We have $H_2(G)=\Z_3$, so $\PP$ is efficient. Moreover $G$ is a Bing group.
\begin{proof}
We will need the following GAP \cite{GAP} program, that uses the packages HAP \cite{HAP} and SONATA \cite{SONATA}.
\begin{verbatim}
LoadPackage("HAP");;
LoadPackage("SONATA");;
F:=FreeGroup(2);;
G:= F/[F.1^3, F.1*F.2*F.1^-1*F.2*F.1*F.2^-1*F.1^-1*F.2^-1,
F.1^-1*F.2^-4*F.1^-1*F.2^2*F.1^-1*F.2^-1];;
Order(G);
G:=SmallGroup(IdGroup(G));;
R:=ResolutionFiniteGroup(G,3);;
Homology(TensorWithIntegers(R),2);
Set(List(Endomorphisms(G),
f->Homology(TensorWithIntegers(EquivariantChainMap(R,R,f)),2)));
\end{verbatim}
The program prints the order of $G$, a list with the invariant factors of $H_2(G)$ and a list with the endomorphisms of $H_2(G)$ that are induced by an endomorphism of $G$. The output is:
\begin{verbatim}
243
[ 3 ]
[ [ f1 ] -> [ <identity ...> ], [ f1 ] -> [ f1 ] ]
\end{verbatim}
Therefore $|G|=243$ and $H_2(G)=\Z_3$. Since for every endomorphism $\phi:G\to G$ we have that $H_2(\phi)$ is either the zero map or the identity, $G$ is a Bing group.
\end{proof}
\end{prop}

By Theorem \ref{BingGroups} and Proposition \ref{GIsBing} we have:
\begin{coro}\label{CoroBing1}
The complex $X_\PP$ associated to the presentation $$\PP = \langle x,y \mid x^3,\, xyx^{-1}yxy^{-1}x^{-1}y^{-1},\, x^{-1}y^{-4}x^{-1}y^2x^{-1}y^{-1} \rangle$$ has the fixed point property. Moreover $\chi(X_\PP)=2$.
\end{coro}

\begin{coro} There are compact $2$-dimensional polyhedra with the fixed point property and Euler characteristic equal to any positive integer $n$.
\begin{proof}
For $n=1$ this is immediate. For $n>1$ take a wedge of $n-1$ copies of the space $X_\PP$ of Corollary \ref{CoroBing1}.
\end{proof}
\end{coro}
To prove Theorem \ref{Bing8Example} we will need another efficient Bing group:

\begin{prop}\label{HIsBing}
The group $H$ presented by $\QQ=\langle x , y \mid x^4, y^4, (xy)^2, (x^{-1}y)^2\rangle$ is a finite group of order $16$. We have $H_2(H)=\Z_2\oplus \Z_2$, so $\QQ$ is efficient. Moreover $H$ is a Bing group.
\begin{proof}
As above we will use a GAP program.
\begin{verbatim}
LoadPackage("HAP");;
LoadPackage("SONATA");;
F:=FreeGroup(2);;
H:= F/[F.1^4, F.2^4, (F.1*F.2)^2, (F.1^-1*F.2)^2];;
Order(H);
H:=SmallGroup(IdGroup(H));;
R:=ResolutionFiniteGroup(H,3);;
Homology(TensorWithIntegers(R),2);
Set(List(Endomorphisms(H),
f->Homology(TensorWithIntegers(EquivariantChainMap(R,R,f)),2)));
\end{verbatim}
The program produces the following output:
\begin{verbatim}
16
[ 2, 2 ]
[ [ f1, f2 ] -> [ <identity ...>, <identity ...> ],
[ f1, f2 ] -> [ f1, f2 ],[ f1, f2 ] -> [ f1^-1*f2^-1, f2^-1 ] ]
\end{verbatim}
This proves that $|H|=16$, $H_2(H)=\Z_2\oplus \Z_2$ and that $H$ is a Bing group.
\end{proof}
\end{prop}

We recall the following theorem:

\begin{teo}[Jiang,{\cite[Theorem 7.1]{Jiang}}]\label{JiangFPPInvariante}
In the category of compact connected polyhedra without global separating points, the fixed point property is a homotopy type invariant.

Moreover, if $X\simeq Y$ are compact connected polyhedra such that $Y$ lacks the fixed point property and $X$ does not have global separating points, then $X$ lacks the fixed point property.
\end{teo}

The following shows that the answer to Question \ref{Bing8} is $2$:

\begin{teo}\label{Bing8Example} There is a compact $2$-dimensional polyhedron $Y$ without the fixed point property and such that the polyhedron $X$, obtained from $Y$ by an elementary collapse of dimension $2$, has the fixed point property.
\begin{proof}
Let $\PP$ and $\QQ$ be the presentations of Propositions \ref{GIsBing} and \ref{HIsBing}. By Theorem \ref{BingGroups}, $X_\PP$ and $X_\QQ$ have the fixed point property, so $X=X_\PP\vee X_\QQ$ also has the fixed point property. Since neither $X_\PP$ nor $X_\QQ$ have global separating points, by adding a $2$-simplex, we can turn $X$ into a polyhedron $Y$, without global separating points and such that, by collapsing that $2$-simplex, we obtain $X$. We have $H_2(\pi_1(Y))= H_2( \pi_1(X_\PP) *\pi_1(X_\QQ))=H_2(\pi_1(X_\PP))\oplus H_2(\pi_1(X_\QQ))=\Z_2\oplus \Z_6$ and $\rk( H_2(Y))=3$. By \cite[Proposition 3.3]{BarmakSadofschi} and Theorem \ref{JiangFPPInvariante}, $Y$ does not have the fixed point property.
\end{proof}
\end{teo}

\end{document}